\documentclass[12pt]{article}
\usepackage{amssymb,amsfonts,epic,graphics,pictex}
\begin{document}

\newtheorem{lem}{Lemma}[section]
\newtheorem{prop}{Proposition}
\newtheorem{con}{Construction}[section]
\newtheorem{defi}{Definition}[section]
\newtheorem{coro}{Corollary}[section]
\newcommand{\hf}{\hat{f}}
\newtheorem{fact}{Fact}[section]
\newtheorem{theo}{Theorem}
\newcommand{\Br}{\Poin}
\newcommand{\Cr}{{\bf Cr}}
\newcommand{\dist}{{\bf dist}}
\newcommand{\diam}{\mbox{diam}\, }
\newcommand{\mod}{{\rm mod}\,}
\newcommand{\compose}{\circ}
\newcommand{\dbar}{\bar{\partial}}
\newcommand{\Def}[1]{{{\em #1}}}
\newcommand{\dx}[1]{\frac{\partial #1}{\partial x}}
\newcommand{\dy}[1]{\frac{\partial #1}{\partial y}}
\newcommand{\Res}[2]{{#1}\raisebox{-.4ex}{$\left|\,_{#2}\right.$}}
\newcommand{\sgn}{{\rm sgn}}

\newcommand{\C}{{\bf C}}
\newcommand{\D}{{\bf D}}
\newcommand{\Dm}{{\bf D_-}}
\newcommand{\N}{{\bf N}}
\newcommand{\R}{{\bf R}}
\newcommand{\Z}{{\bf Z}}
\newcommand{\tr}{\mbox{Tr}\,}

\newenvironment{nproof}[1]{\trivlist\item[\hskip \labelsep{\bf Proof{#1}.}]}
{\begin{flushright} $\square$\end{flushright}\endtrivlist}
\newenvironment{proof}{\begin{nproof}{}}{\end{nproof}}

\newenvironment{block}[1]{\trivlist\item[\hskip \labelsep{{#1}.}]}{\endtrivlist}
\newenvironment{definition}{\begin{block}{\bf Definition}}{\end{block}}

\newtheorem{conjec}{Conjecture}
\newtheorem{com}{Comment}
\newtheorem{exa}{Example}
\font\mathfonta=msam10 at 11pt
\font\mathfontb=msbm10 at 11pt
\def\Bbb#1{\mbox{\mathfontb #1}}
\def\lesssim{\mbox{\mathfonta.}}
\def\suppset{\mbox{\mathfonta{c}}}
\def\subbset{\mbox{\mathfonta{b}}}
\def\grtsim{\mbox{\mathfonta\&}}
\def\gtrsim{\mbox{\mathfonta\&}}

\newcommand{\ar}{{\bf area}}
\newcommand{\1}{{\bf 1}}
\newcommand{\Bo}{\Box^{n}_{i}}
\newcommand{\Di}{{\cal D}}
\newcommand{\gd}{{\underline \gamma}}
\newcommand{\gu}{{\underline g }}
\newcommand{\ce}{\mbox{III}}
\newcommand{\be}{\mbox{II}}
\newcommand{\F}{\cal{F_{\eta}}}
\newcommand{\Ci}{\bf{C}}
\newcommand{\ai}{\mbox{I}}
\newcommand{\dupap}{\partial^{+}}
\newcommand{\dm}{\partial^{-}}
\newenvironment{note}{\begin{sc}{\bf Note}}{\end{sc}}
\newenvironment{notes}{\begin{sc}{\bf Notes}\ \par\begin{enumerate}}%
{\end{enumerate}\end{sc}}
\newenvironment{sol}
{{\bf Solution:}\newline}{\begin{flushright}
{\bf QED}\end{flushright}}

\title{
Multipliers of periodic orbits of quadratic polynomials 
and the parameter plane
}

\author{Genadi Levin
\\
\small{Inst.\ of Math., Hebrew University, Jerusalem 91904, Israel}\\
}
\normalsize
\maketitle
\abstract{We prove an extension result for the multiplier
of an attracting periodic orbit of a quadratic 
map as a function of the parameter. 
This has applications to the problem of geometry
of the Mandelbrot and Julia sets. In particular, we
prove that the size of $p/q$-limb of a hyperbolic component
of the Mandelbrot set of period $n$ is $O(4^n/p)$, and
give an explicit condition on internal arguments
under which the Julia set of corresponding
(unique) infinitely renormalizable quadratic polynomial
is not locally connected.
}

\

\section{Introduction}\label{s1}
Douady-Hubbard-Sullivan (DHS) theorem~\cite{DH1},~\cite{CG} states that
the multiplier $\rho$ of an attracting periodic orbit
is a conformal isomorphism from 
a hyperbolic component of the Mandelbrot set onto the unit disk
$\{|\rho|<1\}$,
and it extends homeomorpically to the boundaries.

In Theorem~\ref{univ} 
we prove that  
$\rho$ extends further to an analytic isomorphism from
a region containing the hyperbolic component onto
a simply connected
domain $\tilde\Omega_n$ containing $\{|\rho|\le 1\}\setminus \{1\}$,
such that the domain $\tilde\Omega_n$ is explicitely defined by the
period $n$ of the attracting orbit. 
This follows from the Main Inequality,
see Theorem~\ref{nunu} below,
which in turn is based upon Theorem 2 of~\cite{Le}, 
see formulation in Subsection~\ref{base} below. 
We derive from Theorem~\ref{univ} few consequences
including stated below Theorems~\ref{limbintro}-~\ref{intro}.

Let us be more precise.
A hyperbolic component $W$ of period $n$ is a component of the interior
of the Mandelbrot set $M$, such that, for $c$ in $W$,  
the quadratic map $f_c(z)=z^2+c$ has an attracting
periodic orbit of period $n$. Denote by $\rho_W(c)$
the multiplier of this orbit of $f_c$, $c\in W$.
By DHS theorem, $\rho_W$ performs a homeomorphism
of the closure of $W$ onto the closed unit disk.
A number $t$ is called an internal argument of a point $c\in \partial W$
iff $\rho_W(c)=\exp(2\pi it)$. The point $c$
with $t=0$ is called the root of $W$.
If $c\in \partial W$ 
is not the root and has a rational internal argument $t$,
the connected component of $M\setminus \{c\}$ which is disjoint
with $W$ is not empty and called the $t$-limb $L(W, t)$ of $W$. 

The first consequence of Theorem~\ref{univ}
concerns the size of the limbs. It
strengthens Yoccoz's bound (off the root), see Section~\ref{limbs},
particularly,
part (b) of Proposition~\ref{lc} and Comment~\ref{gen}.                   .
\begin{theo}\label{limbintro}
There exists $A>0$, such that, for 
every hyperbolic component $W$ of period $n$
and every $t=p/q\in [-1/2, 1/2]\setminus \{0\}$,
the diameter of the limb $L(W, t)$ is bounded by:
\begin{equation}\label{allboundintro}
diam L(W, t)\le A \frac{4^n}{p}
\end{equation}
\end{theo}

This bound immediately implies the local connectivity of the Mandelbrot set
at some parameters $c_*$, where $f_{c_*}$ is infinitely renormalizable
with prescribed unbounded combinatorics, see Corollary~\ref{lcm}.

Let's discuss another result.
Douady and Hubbard proposed an inductive construction
to build an infinitely renormalizable quadratic map
with non locally connected Julia set, and such that
the Mandelbrot set is locally connected at this parameter
~\cite{Sor},~\cite{Mil}.
Their construction consists in 
choosing a sequence of internal arguments $t_m$ 
of successive bifurcations step by step,
by continuity in the corresponding dynamical planes
and in the parameter plane (with help of Yoccoz's bound).
Next statement makes explicit  
the two sides of Douady-Hubbard's construction
(parameter and dynamical). 
\begin{theo}\label{intro}
Let $n\ge 1$. Let
$$t_0, t_1,..., t_m,...$$
be any sequence
of rational numbers $t_m=p_m/q_m\in (-1/2, 1/2]$ which satisfy
the following properties:
if we denote $n_0=n$, $n_m=nq_0...q_{m-1}$, $m>0$, then,
for all $m$ large enough,
\begin{equation}\label{assump1}
p_m > (2B_0)  4^{n_m}, \ \ \ \ 2n_m^2 < \frac{p_m^2}{q_m},
\end{equation}
(where the constant $B_0$ is taken from Theorem~\ref{nunu} below),
and also
\begin{equation}\label{assump2}
\sum_{m=1}^\infty |t_m|^{1/q_{m-1}}<\infty.
\end{equation}
Given a hyperbolic component $W$ of the Mandelbrot set
of period $n$,
consider the following sequence of hyperbolic
components $W^m$:  $W_0=W$, and, for $m>0$,
$W^m$ touches the hyperbolic
component $W^{m-1}$ at the point with an internal argument
$t_{m-1}$.
For every $m$, consider the $t_m$-limb $L(W^m, t_m)$ of $W^m$
(it contains $W^{m+1}$).
Then the limbs $L(W^m, t_m)$ shrink 
to a unique point $c_*$,
the Mandelbrot set is locally connected at $c_*$, and
the map $f_{c_*}$ is infinitely renormalizable
with non locally connected Julia set.
\end{theo}
This statement can be reformulated in term
of the combinatorial data
of a fixed infinitely renormalizable map, see Theorem~\ref{infi}.

Let us comment on the inequalities (~\ref{assump1})-(~\ref{assump2})
above.
Conditions (~\ref{assump1}) guarantee
that the corresponding multipliers are local parameters.
Namely, if $\rho_{W^m}(c)$ is the multiplier of the
periodic orbit of $f_c$, which is attracting for $c\in W^m$,
and $\psi_m$ denotes an inverse to
$\log\rho_{W^m}$, then the first inequality
in (~\ref{assump1}) ensures that $\psi_m$
extends as a holomorphic function
to a disk of radius proportional to $n_m/q_m$ around
$2\pi it_m$
while the second inequality in (~\ref{assump1})
implies that this analytic continuation
is in fact injective.
In turn,
condition (~\ref{assump2}) will guarantee that
the bifurcated periodic orbits stay away from the origin. 
It confirms
Milnor's ``guess'' from~\cite{Mil}, p. 21.

\subsection{Notations}
We collect some notations to be used throughout the paper.

$f_c(z)=z^2+c$, $J_c=J(f_c)$ the Julia set,
$D_{\infty}(c)$ the basin of infinity,

$M=\{c\in {\bf C}: J_c \ \ connected\}$ is the Mandelbrot set,

$B_c$ is the Bottcher coordinate at infinity
normalized by $B_c(z)\sim z$ as $z\to \infty$,

$G_c=\lim_{n\to \infty} 2^{-n}\log |f_c^n(z)|$ Green's function 
of $D_{\infty}(c)$, $G(z)\sim \log|z|$ at $\infty$ extended
by $0$ to the whole plane, so that $G_c(z)=\log B_c(z)$ near infinity.

$B(a,r)=\{z: |z-a|<r\}$.

\section{Main inequality}
Let $O=\{b_k\}_{k=1}^n$ be a periodic orbit 
of $f=f_{c_0}$ of exact period $n$.
Its  multiplier is the number
$\rho=(f^n)'(b_k)$. 
If $\rho\not=1$, by the Implicit function theorem,
there exist $n$ holomorphic
functions $b_k(c)$, $k=1,..,n$ defined in a neighborhood
of $c_0$, such that $O(c)=\{b_k(c)\}_{k=1}^n$
is a periodic orbit of $f_c$, and $O(c_0)=O$.
In particular, if 
$\rho(c)=(f_c^n)'(b_k(c))=2^n b_1(c)...b_n(c)$
denotes the multiplier of $O(c)$, it is holomorphic
in $c$ in this neighborhood. 

In what follows we assume that the multiplier $\rho$ is not very
big because we are interested in studying the behavior
of a multiplier not far from the hyperbolic component
where the corresponding periodic orbit is attracting.
So, we will always assume that 
$$|\rho|<e.$$
\begin{theo}\label{nunu}
There exist $\lambda_*$, $B$, $B_0$ as follows.
Let $O(c)$ be a repelling periodic orbit
of $f_c$ of exact period $n$, and the multiplier 
of $O(c)$ is equal to $\rho$. Then
the following inequality holds:
\begin{equation}\label{main}
|n\rho-n+2c\frac{\rho'(c)}{\rho(c)}|\le K_n(c)
\{\log |\rho(c)|+\frac{|\rho'(c)|}{|\rho(c)|}
\frac{1}{\pi}{\ar}(\{z: 0<G_c(z)<2^{-n}\log\lambda_*\})\}
\end{equation}
Here 
$$K_n(c)=\frac{2}{\log \lambda_*}
\max\{ |(f_c^n)'(z)|: G_c(f_c^n(z))=2\log\lambda_*\},$$
and we have:
\begin{equation}\label{K}
K_n(c)\le B (2\beta)^n
\end{equation}
where $\beta$ is the unique positive solution of the equation
$\beta^2-|c|=\beta$,
and also
$$K_n(c)\le B_0 4^n.$$
\end{theo}

\begin{com}\label{level}
Let us notice for the future use few bounds related to the number $G_c(0)$.

(a) First, the following fact was established in~\cite{EL}:
for the multiplier $\rho$ of every repelling periodic
orbit of $f_c$ of period $n$, we have
$$G_c(0)\le \frac{1}{n}\log |\rho|.$$
In particular, since we assume $|\rho|<e$, then
$c$ belongs to a neighborhood of $M$ where
$G_c(0)\le n^{-1}\le 1$.

(b) Second, the transfinite diameter of the set
$\{z: 0\le G_c(z)<2^{-n}\log\lambda\}$ ($\lambda>1$) is equal to 
$\lambda^{2^{-n}}$, and, hence, by a theorem of Polya~\cite{Gol},
$${\ar}(\{z: 0\le G_c(z)<2^{-n}\log\lambda\})\le \pi \lambda^{2^{-n+1}}.$$
\end{com}

Consider two particular cases of the inequality (~\ref{main}).
One case corresponds to passing to a limit as $|\rho|\to 1$, 
and in another one put $\rho'=0$. Then we get
the following corollary.
\begin{coro}\label{maincoro}

(A) There is $C_0$, such that, if, for some $c_0$, 
$|\rho(c_0)|=1$, then  
\begin{equation}\label{A}
|\rho'(c_0)|\ge C_0\frac{n|\rho(c_0)-1|}{K_n(c_0)}
\end{equation}

(B) For every $n$-periodic orbit $O$ of $f_c$
with the multiplier $\rho$,
if $1<|\rho|<e$ and
\begin{equation}\label{nana}
n|\rho-1| > K_n(c)\log |\rho|,
\end{equation}
then $\rho'\not=0$.
\end{coro}
\begin{com}
(A) implies, of course, that $\rho'\not=0$ on the boundary
of a hyperbolic component where $|\rho|=1$, and it was known
after~\cite{DH2}.
\end{com}

\section{Proof of Theorem~\ref{nunu}}
\subsection{Derivative of multiplier and Ruelle operator}\label{base}
The proof is based on the following result, see Theorem 2 of~\cite{Le}.
Consider the Ruelle transfer operator 
$$Tg(z)=\sum_{w:f_c(w)=z}\frac{g(w)}{(f_c'(w))^2}$$
Given a periodic orbit $O=\{b_k\}_{k=1}^n$ of $f$ of exact
period $n$ and with the multiplier $\rho\not=0,1$, associate to $O$ a function
\begin{equation}\label{A-func}
A(z)=A_O(z)=\sum_{k=1}^n \frac{1}{(z-b_k)^2}+\frac{1}{\rho(1-\rho)}
\sum_{k=1}^n \frac{(f^n)''(b_k)}{z-b_k}
\end{equation}
where $(f^n)''(z)$ is the second derivative of the $n$-iterate
of $f$ w.r.t. $z$.
Then
\begin{equation}\label{ruelle}
A(z)=(TA)(z)+2\frac{\rho'}{\rho}(T1)(z).
\end{equation}

\begin{com}
In~\cite{Le} we use this to give another, {\it local}, proof
of DHS theorem.
\end{com}

\subsection{Two identities}
Denote
$$\gamma_k=\frac{(f^n)''(b_k)}{\rho(1-\rho)}.$$
Let's find the coefficients $C_1$ and $C_2$
in the expansion of the function $A$ at infinity:
$$A(z)=C_1 z^{-1}+C_2 z^{-2} + O(|z|^{-3}).$$
Since
$$A(z)=\sum_{k=1}^n 
(\frac{1}{(z-b_k)^2}+\frac{\gamma_k}{z-b_k}),$$
then
$$C_1=\sum_{k=1}^n\gamma_k, \ \ \ C_2=n+\sum_{k=1}^n b_k\gamma_k$$
On the other hand, comparing the expansions at $\infty$
in (~\ref{ruelle}), one can see that
$$C_1=\frac{\rho'}{\rho},  \  \  \  C_2=2c\frac{\rho'}{\rho}.$$
Thus we have
\begin{equation}\label{iden}
\sum_{k=1}^n\gamma_k,=\frac{\rho'}{\rho},  \ \ \
\sum_{k=1}^n b_k\gamma_k =-n+2c\frac{\rho'}{\rho}
\end{equation}

\subsection{Estimate from above}
Fix $c$ and denote $f=f_c$, $G=G_c$ etc. 
Given $\lambda>1$, denote $U_\lambda=\{z: 0<G(z)<\log \lambda\}$
and $C_\lambda=U_\lambda\setminus f^{-1}(U_{\lambda})$.
Introduce also a number associated to the cycle $O$
as follows: 
\begin{equation}\label{M}
M=\liminf_{\delta\to 0}\int_{B(b,\delta)\setminus f_O^{-n}(B(b,\delta))}
\frac{d\sigma_z}{|z-b|^2}
\end{equation}
where $b$ is a point of the orbit, and
$d\sigma_z$ is the area element in the $z$-plane, 
$b$ is any point in $O$,
and $f^{-n}_O$ is a branch defined near $b$
that fixes $b$.
It is easy to see that 
$$M\le 2\pi\log|\rho|.$$

\begin{com} 
One can show 
(though we'll not use this fact below) that $M=2\pi\alpha\log|\rho|$ where
$\alpha$ is the the density of $D_\infty$ at the point $b$
in the logarithmic metric~\cite{Le1}.
\end{com} 

\begin{lem}\label{1}
For every $\lambda>1$,
\begin{equation}\label{basic}
\int_{C_\lambda}|A(z)|d\sigma_z\le M+2\frac{|\rho'|}{|\rho|}
{\ar}(\{z: 0<G(z)<\frac{1}{2}\log\lambda\})
\end{equation}
\end{lem}
\begin{proof}
Denote by $f_O^{-i}$ a local branch sending
a point $b$ of the orbit $)$ 
to $f^{n-i}(b)$, $i=1,2,...,n$. For $\delta>0$, let
$V_\delta=U_\lambda\setminus \cup_{i=0}^{n-1} f_O^{-i}(B(b,\delta))$.
Then $A$ is integrable in $V_\delta$, and one can write
$$\int_{V_\delta}|A|d\sigma= 
\int_{V_\delta}|TA(z)+2\frac{\rho'}{\rho} T1(z)|d\sigma\le
\int_{f^{-1}(V_\delta)}|A|d\sigma+2\frac{|\rho'|}{|\rho|}
{\ar}(f^{-1}(V_\delta)).$$
Note that 
$V_\delta= f^{-1}(V_\delta)\cup C_{\lambda}\cup \Delta\setminus 
(B(b,\delta)\setminus f_O^{-n}(B(b,\delta))$
where $\Delta$ is an open set, which shrinks
to $-O=(-b_1,...,-b_n)$ as $\delta\to 0$.
Therefore, 
$$\int_{C_{\lambda}}|A|d\sigma\le 
\liminf_{\delta\to 0}\int_{B(b,\delta)\setminus 
f_O^{-n}(B(b,\delta))}|A(z)|d\sigma_z+2\frac{|\rho'|}{|\rho|}
{\ar}(U_{\lambda^{1/2}}).$$
We have: $A(z)=(z-b)^{-2}+r(z)$ where $r(z)$ is integrable 
at $b$.
Sending $\delta\to 0$, we get the desired inequality. 
\end{proof}

\subsection{Estimate from below}
We are left with the problem to estimate the integral
\begin{equation}\label{int}
I(\lambda):=\int_{C_\lambda}|A(z)|d\sigma_z
\end{equation}
from below. 

To be able to deal with $I(\lambda)$ when $J$ is disconnected
and $\lambda$ close to $1$, let us extend the Bottcher function
$B$ of $f$ at $\infty$ to the following simply connected 
domain (on the Riemann sphere) $D_\infty^*$: it is 
obtained from $D_\infty$ by deleting 
all arcs of external rays (gradient curves of $G$) starting
from $0$ and
all its preimages up to the Julia set.
Denote the extended Bottcher function again by $B$.
It maps $D_\infty^*$ onto a domain $A$ whose 
boundary is a 'hedgehog', see~\cite{LeSo},~\cite{Leyo}.
Let $\phi: A\to D_\infty^*$ be the inverse map.
Note that $\phi'$ has a singularity at the
vertex $p$ of every needle (i.e. the image by $B$ of a 
critical point of $G_c$), but it is of the type
$|w-p|^{-1/2}$, i.e.
integrable against the area.

Let $\lambda>1$. One writes:
$$I(\lambda)=\int_{C_\lambda}|A(z)|d\sigma_z=
\int_{\lambda<|w|<\lambda^2} |A(\phi(w))||\phi'(w)|^2 d\sigma_w.
$$
Let us assume that $f^n(\phi(w))\not=0$ 
for every $\lambda\le |w|\le \lambda^2$, 
in other words, that $\lambda^{2^n}>G(0)$.
By Fubini theorem,
we can continue as follows:
$$\int_{\lambda<|w|<\lambda^2} |A(\phi(w))||\phi'(w)|^2 d\sigma_w\ge
\int_{\lambda}^{\lambda^2}\min_{|w|=r}\frac{|\phi'(w)|}{|f^n(\phi(w))|}dr
\int_{\{z: G(z)=\log r\}}|A(z)f^n(z)||dz|.$$
Now,
$$\int_{\{z: G(z)=\log r\}}|A(z)f^n(z)||dz|\ge 
|\int_{\{z: G(z)=\log r\}}A(z)f^n(z) dz|=2\pi |R|,$$
where
$$R=\frac{1}{2\pi i}\int_{\{z: G(z)=\log r\}}A(z)f^n(z) dz.$$
By the definition of $A$ and the Cauchy formula,
$$R=\sum_{k=1}^n ((f^n)'(b_k)+\gamma_k f^n(b_k))=
n\rho-n+2c\frac{\rho'}{\rho},$$
where we use the identity (~\ref{iden}).
Thus, for every $\lambda>1$, such that $\lambda^{2^n}>G(0)$, 
\begin{equation}\label{est}
I(\lambda)\ge 2\pi |R|\int_{\lambda}^{\lambda^2}
\min_{|w|=r}\frac{|\phi'(w)|}{|f^n(\phi(w))|}dr.
\end{equation}
Now, from $f\phi(w)=\phi(w^2)$ we conclude that
$$\frac{\phi'(w)}{f^n(\phi(w))}=
\frac{2^n w^{2^n-1}\phi'(w^{2^n})}{\phi(w^{2^n})(f^n)'(\phi(w))}.$$
There is a choice for $\lambda>1$.
Take $\lambda=\lambda_n$ where $\lambda_n^{2^n}=\lambda_*$ and
$\lambda_*$ is fixed. 

Let us estimate $|c|$.
By Comment~\ref{level}, (a),
$G_c(0)\le n^{-1}$. 
If $J_c$ is connected, then $|c|\le 2$.
Otherwise
the function $\Psi_c$ which is inverse to the Bottcher
coordinate $B_c$ of $f_c$ extends in a univalent fashion to the disk
$\{|w|>\exp(G_c(0))\}$. Besides, $c=\Psi(w_c)$, for some
$|w_c|=\exp (2 G_c(0))$, and $\Psi$ is odd.
Therefore, by a classical
distortion theorem, see e.g.~\cite{Gol},
\begin{equation}\label{c}
|c|=|\Psi(w_c)|\le 2 \exp(2 G_c(0))\le 2 \exp(\frac{2}{n}), n>0.
\end{equation}
Thus $c$ belongs to a bounded neighborhood of the Mandelbrot set.
Hence, we can fix $\lambda_*$ in such
a way
that for all $c$ in the fixed neighborhood,
$|\phi'(w)/\phi(w)|>1/(2|w|)$  
for all $\lambda_*<|w|<\lambda_*^2$.
We define
$$K_n(c)=\frac{2\max\{ |(f_c^n)'(z)|: G_c(z)=2^{-n+1}\log\lambda_*\}}
{\log \lambda_*}.$$
The inequality (~\ref{main}) follows.

\subsection{The constant}
It remains to estimate $K_n(c)$.
Let us show that there is $B$ such that
\begin{equation}\label{KK}
K_n(c)\le B (2\beta)^n
\end{equation}
where $\beta$ is a unique positive solution of the equation
$\beta^2-|c|=\beta$. 
(Note that $\beta>1$ because $c\not=0$.)
Indeed, if, for some $\delta>0$,
$|z|=\beta+\delta$, then $|f_c(z)|>\beta + 2\beta\delta$.
Hence, if $|f_c^n(z)|=R$ where $R>1+\beta$, then it follows
$|f_c^k(z)|-\beta< (R-\beta)(2\beta)^{k-n}$.
Let us choose and fix $R$, such that $R>1+\beta$ 
and $\{z: G_c(z)=2\log\lambda_*\}\subset
B(0, R)$ for every $c$ from the neighborhood of $M$.
Then we get (~\ref{KK}) with
$$B=\frac{2}{\log\lambda_*}
\Pi_{k=1}^\infty (1+\frac{R}{2^k}).$$ 

On the other hand, as we know,
$$|c|\le 2 \exp(\frac{2}{n})
\le 2(1+\frac{A_1}{n}),$$
for some $A_1$.
Therefore,
$$\beta\le 2+\frac{A_2}{n}$$
and, hence,
$$K_n(c)\le B 4^n (1+\frac{A_2}{2n})^n\le 4^n B_0,$$
for some $A_2$ and $B_0$.

\section{Analytic extension of the multiplier}
\subsection{A set of domains}
Given $C>1$, consider an open set of points $\Omega$ 
in the punctured $\rho$-plane defined by the inequality
\begin{equation}\label{domeq}
|\rho-1| > C\log|\rho|
\end{equation}
It obviously contains the set $D_*=\{\rho: 0<|\rho|\le 1, \rho\not=1\}$
and is disjoint with an interval $1<\rho<1+\epsilon$.
Denote by $\Omega_C$ a connected component of $\Omega$
which contains the set $D_*$ completed by $0$. 
Denote also by $\Omega_C^{\log}$ the set of points 
$$L=\log \rho=x+iy, \ \ \rho\in \Omega_C, \ \ |y|\le \pi$$
Note that $\Omega_C^{\log}$ is symmetric w.r.t. the real axis.
\begin{lem}\label{dom}
$\Omega_C$ is simply-connected.
More precisely, the intersection of $\Omega_C^{\log}$ with
any vertical line with $x=x_0>0$  is either empty
or equal to two (mirror symmetric) intervals.
If $C>2$, then $x< 2/(C-2)$ for all $L=x+iy$ in $\Omega_C^{\log}$.
If $C$ is large enough, 
$\Omega_C^{\log}$ contains two (mirror symmetric) domains
bounded by the lines $y=\pm (C/2)x$ ($x>0$) and $y=\pm \pi$
\end{lem}
\begin{proof}
In the log-coordinate
$L=\log \rho=x+iy$, $|y|\le \pi$, the condition
(~\ref{domeq}) with $x>0$ is equivalent to
\begin{equation}\label{logdom}
\sin^2(\frac{y}{2})>\frac{C^2x^2-(\exp(x)-1)^2}{4\exp(x)}
\end{equation}
Given $x$, the set of $y$ which satisfies the latter inequality
is either empty or a union
$[-\pi,-y)\cup (y,\pi]$ with some $y>0$.
We have for $x+iy\in \Omega_C^{\log}$: 
$1\ge (C^2x^2-(\exp(x)-1)^2)/4\exp(x)$, and since $x>0$, then
$Cx<\exp(x)+1$. If $C>2$, it implies  $x< 2/(C-2)$. 
If $C$ large, then by this $x>0$ must be small, 
and the boundary curve of $\Omega_C^{\log}$
can be written in the form
$\sin(y/2)=\pm (1/2)(C^2-1)^{1/2} x (1-x/2 + x r(x,C))$ where
$r(x,C)=o(1/C^2)$ uniformly in $x$. The statement follows. 
\end{proof}

\subsection{Hyperbolic components and related definitions}
Let $W$ be a component of the interior of $M$.
It is called a hyperbolic component if 
$f_c$, $c\in W$, has an attracting periodic orbit $O(c)$.
Let us call $W$ $n$-hyperbolic, if the
exact period of the latter orbit is $n$.

Denote by $\rho_W(c)$ the 
multiplier of $O(c)$.
By the DHS theorem, $\rho_W$
is a analytic isomorphism of $W$ onto the unit disk,
and it extends homeomorphically to the boundaries.
Given a number $t\in (-1/2, 1/2]$, denote by
$c(W, t)$ the unique point in $\partial W$ with the
{\it internal argument} $t$, i.e.
$\rho_W$ at this points is equal to $\exp(2\pi t)$. 
Root of $W$ is the point $c_W=c(W, 0)$ with the internal argument zero.
$W$ is called primitive iff its root is not a point of other
hyperbolic component.

If $t=p/q$
is a rational number, we will always assume that $p, q$ are 
co-primes.
For any rational $t\not=0$, denote
by $L(W, t)$ the connected component 
of $M\setminus \{c(W, t)\}$ which is disjoint with $W$. 
It is called the $t$-limb of $W$.
Denote also by $W(t)$ a $nq$-hyperbolic component
with the root point $c(W, t)$; it
touches $W$ at this point. The limb $L(W, t)$ contains $W(t)$.
The root $c_W$ of $W$ is the landing point of precisely two
external rays of the Mandelbrot set~\cite{DH2}.
In what follows, 
it will be important the notion
of the wake of a hyperbolic component $W$~\cite{H}: 
it is the only component $W^*$ of the plane cut by two external rays to
the root of $W$ that contains $W$.
We'll use also the following fact
(see e.g.~\cite{Leyo}, Theorem 7.2): 
\begin{prop}\label{wake}
The points of periodic orbit
$O(c)$ as well as its multiplier
$\rho_W$ extend as analytic functions to the wake $W^*$.
Moreover, $|\rho_W|>1$ in $W^*\setminus \overline W$.
\end{prop}
The map $\rho_W$ from $W$
onto the unit disk
$c\mapsto \rho_W(c)$ has an inverse, which we denote by
$c=\psi_W (\rho)$. It is defined so far in the unit disk.
\subsection{Extension of the multiplier}
Introduce
$$\Omega_n=\Omega_{n^{-1} 4^n B_0}=
\{\rho: |\rho-1|>\frac{4^n B_0}{n}\log|\rho|\}$$
and
$$\Omega_n^{\log}=\Omega_{n^{-1} 4^n B_0}^{\log}=\{L=x+iy: 
\exp(L)\in \Omega_n, |y|\le \pi\}.$$ 
\begin{theo}\label{univ}
 
(a) The function $\psi=\psi_W$ extends to a holomorphic function
in the domain $\Omega_n$.

(b) The function $\psi$ is univalent
in a subset $\tilde\Omega_n$ of $\Omega_n$
defined by its log-projection
$\tilde\Omega_n^{\log}=\{\log\rho: \rho\in \tilde\Omega_n\}$ 
as follows:
$$\tilde\Omega_n^{\log}=\Omega_n^{\log}\setminus
\{L: |L-R_n|<R_n\}$$ 
where $R_n$ depends on $n$ only and has an asymptotics
$$R_n=(2+O(2^{-n}))n\log 2$$ 
as $n\to \infty$. 

Finally, the image of $\tilde\Omega_n$ by $\psi$ is contained
in the wake $W^*$.
\end{theo}

\begin{com} 
The disk $\{L: |L-R_n|<R_n\}$ cuts off from $\Omega_n^{\log}$ 
an asymptotically negligible portion:
the deleted part is contained in the
disk $\{|L|<r_n\}$, where $r_n\sim (4\log 2/B_0) (n^2/4^n)$. 
\end{com}  

\begin{proof}
Let $\rho$ be a multiplier of some
repelling orbit of period $n$ for some $f_c$.
By Corollary~\ref{maincoro}, part (B), 
if $\rho\in \Omega_n$,
then $\rho'\not=0$.

(a) Since $\Omega_n$ is simply connected,
it is enough to show that $\psi$
has an analytic extension along any curve
in $\Omega_n$ (which starts at $0$).

Firstly, $c$ is an algebraic function of $\rho$.
Indeed, consider two functions
$Q(c,b)=f_c^n(b)-b$ and $P(c,b,\rho)=b f_c(b)...f_c^{n-1}(b)-\rho/2^n$.
They
are polynomials in $c,b$ and $\rho$, which are of degrees
$2^n$ and $2^n-1$ respectively w.r.t. $b$ and with 
leading coefficients $1$. Hence, resultant
$R(c,\rho)$ of $Q$ and $P$ w.r.t. $b$ is a polynomial
such that $R(c,\rho)=0$ if and only if $Q$ and $P$ have a joint
root $b$, which means that $b$ is a fixed point of
$f_c^n$ and $\rho=(f_c^n)'(b)$.

Now assume that $\psi$ does not have an analytic continuation
along a (simple) curve $\beta$ in $\Omega_n$ starting in $0$.
By the above, it means that $\beta$ contains
a singular point $\rho_0$ of the algebraic function $c$,
that is, when we make a small loop around $\rho_0$ then we get
a different value of $c$. Let $\rho_0$ be the first such point
when we move from $0$. 
Denote $c_0=\psi(\rho_0)$ the limit value
of $\psi$ when $\rho$ approaches $\rho_0$ along $\beta$.
Since $c$ is algebraic, 
$c$ and $\rho$ can be written in a form:
$$c-c_0=u^j, \rho-\rho_0=(\phi(u))^k$$
where the integers $j,k>0$ and  $\phi$ is holomorphic near $u=0$
such that $\phi(0)=0, \phi'(0)\not=0$.
Note that $\rho_0\not=1$, hence $\rho$ is holomorphic
in $c$ in a neighborhood of $c_0$.
By the same reason, the period of the corresponding periodic orbit
of $f_{c_0}$ (=limit of the orbit of $f_{\psi(\rho)}$
when $\rho$ tends to $\rho_0$ along $\beta$)
is exactly $n$. 
$\rho\in \Omega_n$, hence, $\rho'(c_0)\not=0$.
But $dc/d\rho\sim (j/k) u^{j-k}/\phi'(0)$, which 
implies $j=k$.
Since $\rho$ is holomorphic in $c$ near $c_0$, we get
$$(\phi(u))^k=\tilde\phi(u^k)$$
for another holomorphic near $u=0$
function $\tilde\phi$. Here $\tilde\phi'(0)\not=0$ because
$\phi'(0)\not=0$. Thus in a new local coordinate
$\tilde u=u^k$ we have $c-c_0=\tilde u$, $\rho-\rho_0=\tilde \phi(\tilde u)$.
It follows $\psi$ can be extended through $\rho_0$,
a contradiction.

(b) $\psi$ is univalent in some domain if
it takes values in a simply connected domain 
in the $c$-plane where the function $\rho(c)$
(local inverse to $\psi$)
is well-defined. Let us choose as such a domain
in the $c$-plane the wake of $W$. 
The choice is correct
because as we know 
the function 
$\rho$ extends to a holomorphic function from $W$ to
its wake $W^*$.
Now it is enough to show that, for any $L=\log\rho\in \tilde\Omega_n^{\log}$,
the value $\psi(\rho)$ cannot belong to the boundary of
$W^*$. The latter consists of two external rays in the $c$-plane
to the root of $W$ (including the root itself). 
Assume it is not the case. Then we find
a curve $l$ in $\tilde{\Omega_n}$, which starts at $0$
and ends at some $\rho_0$, such that
$\psi(\rho)\in W^*$ for all $\rho\in l\setminus \{\rho_0\}$
while 
$c_0=\psi(\rho_0)$,
for $L_0=\log\rho_0\in \Omega_n^{\log}$ does belong
to such a ray. Consider a continuation
of the $n$-periodic orbit $O(c)$ of $f_c$ 
along the curve $l$
which is attracting for $c\in W$. Then 
the rotation number 
of the periodic orbit $O(c_0)$ of $f_{c_0}$ is zero
(see e.g. Remark 7.2 of~\cite{Leyo}). Now, by
Theorem 5.1 of ~\cite{LeSo}, 
\begin{equation}\label{cut}
\frac{|L_0|^2}{Re(L_0)}\le 
\frac{2\pi n\log 2}{\pi/2-\arctan[(2^n-1)a_0/\pi]}
\end{equation}
where $a_0=G_{c_0}(0)$.
Let's estimate $a_0$ from above.
According to~\cite{EL}, Theorem 1.6, 
$Re(L_0)=\log|\rho_0|\ge n a_0$.
By Lemma~\ref{dom}, $Re(L_0)<2/(C_n-2)$, where 
$C_n=4^n B_0/n$. Hence,
$a_0<2/(4^n B_0-2n)$.
It allows us to define
$$R_n=\frac{\pi n\log 2}{\pi/2-\arctan[\frac{2(2^n-1)}{\pi(4^n B_0-2n)}]}.$$
Note that $R_n=(2+O(2^{-n}))n\log 2$ as $n\to \infty$.
Now, (~\ref{cut}) implies that 
$|L_0-R_n|<R_n$, i.e. $\rho_0$ belongs to the part of
$\Omega_n$ that we delete, a contradiction. 
\end{proof}

\section{Limbs}\label{limbs}
\subsection{Yoccoz's circles}\label{yoc}
Let $W$ be an $n$-hyperbolic component.
As we know, the multiplier $\rho_W$ is defined and analytic throughout
the wake $W^*$.
Let us formulate a result due to Yoccoz, which is
a basic point for further considerations.
For every $t=p/q\not=0$, consider the limb $L(W,t)$.
Then, for every $c\in L(W,t)$,
a branch of $\log\rho_W(c)$ is
contained in the following round disk (Yoccoz's circle):
\begin{equation}\label{yoccircle}
Y_n(t)=\{L: |L-(2\pi it+\frac{n\log 2}{q})|<\frac{n\log 2}{q}\},
\end{equation}
see~\cite{H},~\cite{Pe} as well as  
~\cite{Po},~\cite{Le},~\cite{EL1},~\cite{EL},~\cite{Leyo}. 


As a well known corollary~\cite{H}, we have:
\begin{prop}\label{lc}

(a) The intersection of the wake $W^*$ of $W$ (completed by the root)
with $M$
is equal to the union of $\overline W$ and its limbs.

(b) For every hyperbolic component $W$ there exists
a constant $C_W$ depending on $W$, such that, if $W$
is not primitive, then 
$$diam L(W, p/q)\le \frac{C_W}{q}$$
and if $W$ is primitive, then
$$diam L(W, p/q)\le \frac{C_W}{q^{1/2}}.$$
\end{prop}

\subsection{Condition on the internal argument}
\begin{defi}\label{deep}
Given $n$, let us call a rational numbers $t=p/q\in [-1/2, 1/2]$
$n$-deep, iff,
for every $n'\le n$ and every $n'$-hyperbolic component $W$,
there is a disk $B(2\pi it, d)$ with $d<\pi$, such that the following holds:

(a) the image $\exp(B(2\pi it, d/2))$ 
of $B(2\pi it, d/2)$ by the exponential map covers the limb $L(W, t)$,

(b) the inverse $\psi_W$ to the map $\rho_W$ extends to a univalent function
defined in the union of the
unit disk and the domain $\exp(B(2\pi it, d))$.
\begin{com}\label{conj11}
By (\ref{yoccircle}),
$$d<\frac{4n\log 2}{q}.$$
Conjecturely,
$d$ should be comparable to $n/q^2$. For simplicity, we will 
always assume that
$$d>\frac{n}{q^4}.$$
\end{com}
\end{defi}
\begin{com}\label{actual}
By the proof of Theorem~\ref{univ}, 
if $\psi_W$ extends just {\it holomorphically}
into $\exp(B(2\pi it, d))$ and
maps it into the wake $W^*$, then $t$ is $n$-deep.
In particular, if 
$B(2\pi it, (4n\log 2)/q)\subset \tilde\Omega_n^{\log}$, then
$t$ is $n$-deep.
\end{com}

\begin{prop}\label{propdeep}

(1) For every fixed $n$,

(1a) $n$-deep rationals are dense
in $(-1/2, 1/2)$; 

(1b) the set of all $t$, which are not $n$-deep has
a unique concentration point $0$.

(2) there exists $n_0$, such that, for every $n>n_0$,
the point $t=p/q\in (-1/2, 1/2)$ is $n$-deep
if the following two inequalities hold:
\begin{equation}\label{first}
p> 4^n (2B_0), 
\end{equation}
\begin{equation}\label{second}
\frac{p^2}{q}>2 n^2. 
\end{equation}
\end{prop}
\begin{proof}
(1) follows from Comment~\ref{actual}. As for (2), by 
the same Comment
it is enough to check that the
inequalities (~\ref{first})-(~\ref{second}) guarantee
that the disk of radius $(4n\log 2)/q$ centered
at the point $2\pi it$ is contained in the domain $\tilde\Omega_n^{\log}$.
Denote $\theta=2\pi t$ and
consider the disk $B(i\theta, r)$. 
Then for $C>0$ and $L=i\theta+w$ where $|w|<r$, $r<1$,
one can write:
\begin{equation}\label{eq}
|\exp(L)-1|-C Re(L)\ge |\exp(-i\theta)-1|-|exp(w)-1|-C|w|\ge
2|\sin(\theta/2)|-(C+2)r.
\end{equation}
That is, if $r<2|\sin(\theta/2)|/(C+2)$, then
the disk $B(i\theta, r)$ lies inside of $\Omega_C^{\log}$.
Note that $|\sin(\theta)|\ge 2\theta/\pi$.
This shows  that if we put here $C=C_n=4^n B_0/n$ and 
$r=r_n=(4n\log 2)/q$, then
the inequality (~\ref{first}) ensures for $n$ big that
$B(i\theta, r_n)$ lies in $\Omega_n^{\log}$.

It is also easy to check that the inequality (~\ref{second})
means that $B(i\theta, r_n)$ is disjoint with the disk
$B(R_n, R_n)$, for $n$ big, so that both inequalities imply that
$B(i\theta, r_n)$ is contained in $\tilde\Omega_n^{\log}$.
\end{proof}
\subsection{Uniform bound on the size of some limbs}

\begin{theo}\label{limb}
(1) Let $W$ be an $n$-hyperbolic component, and let
$c\in \partial W$ have an internal argument
$t$, such that $t$ is $n$-deep. Then
there is a topological disk $B(c)$, such that:

$B(c)$ does not contain the root $c_W$ of $W$ and it
is ``roughly'' 
a round disk around the point $c$: for some $r$,
$$B(c, r/4)\subset B(c)\subset B(c, 4 r).$$
The function $\log\rho_W$ 
is univalent in $B(c)$ and maps it
onto a disk $B(2\pi it, d/2)$.
Moreover, $\log\rho_W$ extends univalently to 
a topological disk containing $B(c)$,
and maps this bigger domain onto $B(2\pi it, d)$.

The limb $L(W, t)$ is contained in $B(c)$. 

(2) There exists $A>0$, such that, for 
every $n$-hyperbolic component $W$ and every $t=p/q\in [-1/2, 1/2]$,
the diameter of the limb $L(W, t)$ is bounded by:
\begin{equation}\label{allbound}
diam L(W, t)\le A \frac{4^n}{p}= A \frac{4^n}{t}\frac{1}{q}
\end{equation}
\end{theo}
\begin{com}\label{gen}
Part (2) is Theorem~\ref{limbintro} announced in the Introduction.
It strengthens the bound of part (b)
of Proposition~\ref{lc} (off the root, i.e. if $t=p/q$ is outside
of a neighborhood of zero).
\end{com}
\begin{proof}
(1) Following the Definition~\ref{deep},
one can take $B(c)=\psi_W(\exp(B(2\pi it, d/2)))$.
We use here and later on the following classical distortion bounds
for univalent maps, see e.g.~\cite{Gol}:
if $g$ is univalent in a disk $B(0, R)$ and $r<R$,
then
\begin{equation}\label{distor}
B(g(0), \alpha(r/R)^{-1} r|g'(0)|)\subset g(B(0, R))\subset
B(g(0), \alpha(r/R) r|g'(0)|)
\end{equation}
where $\alpha(x)=(1-x)^{-2}$. 

Now (1) is an immediate consequence 
of the definition~\ref{deep} 
and the latter bound (with $r=R/2$).

(2) Let us prove the second part of the Theorem.
It is enough to show that there exist $M, N$, such that
if $p> M 4^n$, then $diam L(W, t)\le N 4^n/p$ for all $n, t$.
Fix $n$ and $t=p/q$.
In the course of the proof $A_i$ will denote different constants
from a finite collection of numbers.
For $r_n=2n\log 2/q$ and $C_n=n^{-1} 4^n B_0$ one can write:
$$\frac{r_n(C_n+2)}{|\sin(\pi p/q)|}\le A_1\frac{4^n}{p}$$
Hence, if $A_1 4^n/p< 1/4$, then by (~\ref{eq})
the disk $B(2\pi i p/q, 2r_n)$ is contained in $\Omega_n^{\log}$.
Hence, so is the disk
$Y:=\{|L-(i\theta+n\log 2/q)|<n\log 2/q\}$ which is contained in 
$B(2\pi i p/q, 2r_n)$.
By Theorem~\ref{univ}, $\psi$ extends to 
a holomorphic function in $B(2\pi i p/q, 2r_n)$.
For every $c$ in the limb $L(W,t)$ there
is $L\in Y$ such that $c=\psi(L)$.
Therefore, the distance between the
root of the limb and $c$
is bounded from above by the integral
$J:=\int_{\i\theta}^L |\psi'(w)| |dw|$ 
for some $L$ with $|L|=r_n$. On the other hand, by
Theorem~\ref{nunu},
$$|\psi'(w)|\le \frac{|2c n^{-1}|+\pi^{-1}
{\ar}(\{z: G_c(z)<a 2^{-n}\})C_n}
{|\exp(w)-1|-C_n Re(w)}$$
Here area 
${\ar}(\{z: G_c(z)<a 2^{-n}\})$ 
is smaller than $\pi(1+o(1))$ as $n\to \infty$ by
Comment~\ref{level} (b).
Using (~\ref{eq}) we can write
$$J\le 
(1+o(1))C_n r_n\int_{0}^1 \frac{d\tau}{2|\sin(\theta/2)|-r_n(C_n+2)\tau}$$
and hence
$$J\le (1+o(1))C_n \log\frac{2|\sin(\theta/2)|}
{2|\sin(\theta/2)|-r_n(C_n+2)}.$$
Since $r_n(C_n+2)/|\sin(\theta/2)|<1/4$ and 
$\log (1-x)^{-1}\le 2x$ if $0<x<1/4$, we can proceed: 
$$J\le (1+o(1))C_n \frac{r_n}{|\sin(\theta/2)|}\le A_2 \frac{4^n}{p}$$
for all $n$.
Thus it is enough to put
$M=4 A_1$ and $N=A_2$.
\end{proof}
We have the following corollary.
\begin{coro}\label{lcm}
For any sequence
$$t_0, t_1,..., t_m,...$$
of rational numbers $t_m=p_m/q_m\in [-1/2, 1/2]$ the following holds.
Let $W$ by a $n$-hyperbolic component.
Let $W^0=W$, $W^{m}=W^{m-1}(t_{m-1})$, $m=1,2,...$
(i.e. the hyperbolic component $W^m$ touches the hyperbolic
component $W^{m-1}$ at a point with the internal argument
$t_{m-1}$), so that $W^m$ is $n_m$-hyperbolic component
where the periods $n_0=n$, $n_m=nq_0...q_{m-1}$.
Assume that
\begin{equation}\label{assump}
\lim_{m\to \infty}\frac{p_m}{4^{n_m}}=\infty.
\end{equation}
Then the limbs $L(W^m, t_m)$ shrink 
to a unique point $c_*$.
The Mandelbrot set is locally connected at $c_*$.
\end{coro}
\begin{proof}
By the previous result, $diam L(W^m, t_m)\to 0$.
Hence, $L(W^m, t_m)$ shrink to a point $c_*$.
The local connectivity follows because these limbs 
form a shrinking sequence of connected neighborhoods of $c_*$.   
\end{proof}

\section{Selecting internal arguments}
\subsection{Bifurcations}
Throughout this subsection, we consider the following situation.
Let $W$ be a $n$-hyperbolic component,
and let $c_0\in \partial W$ have an internal argument
$t_0=p/q\not=0$. 

Let $O(c)=\{b_j(c)\}_{j=1}^n$ be the periodic orbit
of $f_c$ which is attracting when $c\in W$.
Then all $b_j(c)$ as well as the mulptiplier $\rho(c)$
are holomorphic in $W$ and
extend to holomorphic functions in $c$ near $c_0$
(in fact, in the whole wake of $W$).
As we know the multiplier $\rho(c)$ is injective
near $c_0$. Denote the inverse function
by $\psi(\rho)$, $\psi(\exp(2\pi t_0))=c_0$.
$\psi$ is well defined in $\tilde\Omega_n$, that
includes the unit disk and a neighborhood of the point
$\rho_0=\exp(2\pi it_0)$.

Consider now the wake of the $nq$-hyperbolic component $W(t_0)$
(by definition, it tangents $W$ at the point $c_0$).
By Douady-Hubbard theory~\cite{DH2}, the root $c_0$
of $W(t_0)$ is the landing point of precisely two
external rays of $M$. Denote their arguments by $\tau(c_0)$, $\tilde\tau(c_0)$.
For every $c$ in the wake of $W(t_0)$, we have the
following picture in the dynamical plane
of $f_c$: two external rays of $f_c$ with the arguments
$\tau(c_0)$ and $\tilde\tau(c_0)$ land at one point 
(denote it $b_n(c)$), which
is a point of the repelling periodic orbit $O(c)$.
It implies also that  
$c$ lies in a component of 
the plane minus these two external rays which does not contain the origin.
Indeed, by the formula for
the uniformization map of the exterior 
of $M$~\cite{DH1},
this is true for those $c$ in the wake of $W(t_0)$ which are outside
of the Mandelbrot set. Hence, it must be true throughout
the wake because $c$ cannot cross the external rays
as well as their landing point $b_n(c)$.

We will make use also of a well known formula, see e.g.~\cite{Leyo}:
\begin{equation}\label{digits}
|\tau(c_0)-\tilde\tau(c_0)|=\frac{(\beta-\alpha)(2^n-1)}{2^{nq}-1}
\end{equation}
where $\beta, \alpha\in \{0,1,2,...,2^n-1\}$ are two ``digits''
determined by $W$.

Next statement is combinatorial (cf.~\cite{Che}).
\begin{lem}\label{combin}
Let $c$ be a point of a limb $L(W, t')$ with some
$t'=p'/q'$ and $q'>2$. Assume that $f_c$ has a periodic point
of period $nQ$ with the multiplier $1$. Then
$$Q\ge q'-1$$
\end{lem}
\begin{proof}
Consider the dynamical plane of $f_c$. 
The critical value $c$ belongs to a petal of $a$.
Since $c$ lies in the sector bounded by the two external rays with arguments
$\tau(c')$, $\tilde\tau(c')$ where $c'$ is the root of $L(W, t')$,
then $a$ is in the same sector, too. 
On the other hand, $a$ is a landing point of two external rays 
fixed by $f_n^{nQ}$. Therefore, there must be
$|\tau(c')-\tilde\tau(c')| > (2^{nQ}-1)^{-1}$.
If we apply now the formula (\ref{digits}), it gives
$nQ\ge nq'-2n+1$, that is, $Q\ge q'-1$.   
\end{proof}

Next Lemma describes the bifurcation near the parameter $c_0$, 
cf.~\cite{Che},~\cite{Gu}.
\begin{lem}\label{local} 
There exist a small disk $U$ around the origin,
a neighborhood $V$ of the cycle $O(c_0)$, and $n$ functions
$F_k(s)$, $k=1,...,n$, which are holomorphic
in $U$, such that
$F_k(0)=0$, $F_k'(0)\not=0$ and,
for every $s\in U$ and every $\rho=\rho_0+s^q$, the points 
$$b_k(c_0)+F_k(s\exp(2\pi j/q)), k=1,...,n, j=0,...,q-1$$
are the only fixed points of $f_c^{nq}$ in the
neighborhood $V$ different from $O(c)$, where
$c=\psi(\rho_0+s^q)$ and $\rho$ is the multiplier of $O(c)$.
They form a periodic orbit of $f_c$ of period $nq$.
\end{lem}
\begin{proof}
Introduce
$$p(c,z)=\frac{f_c^{nq}(z)-z}{f_c^n(z)-z}.$$
It is a polynomial in $z$ and $c$.
As we also know, the function $c=\psi(\rho)$
satisfies another polynomial equation $R(c, \rho)=0$,
see the proof of Theorem~\ref{univ}.
Hence, periodic points of period $n$ form
an algebraic function in $\rho$:
they satisfy a polynomial equation
of the form 
\begin{equation}\label{alg}
\tilde p(\rho, z)=0.
\end{equation}
For every $k$, the point $(\rho_0, b_k(c_0))$
with $\rho_0=\exp(2\pi it_0)$ is a singular point
of $\tilde p$. Hence, there exist a local  
uniformizing parameter $s$ and co-prime $i,j$,
such that every solution of (~\ref{alg}) 
in a small enough neighborhood of the singular point
has a form
$$\rho-\rho_0=s^i, z-b_k(c_0)=F_k(s)$$
where
$$F_k(t)=r_k s^j+O(s^{j+1})$$
is a holomorphic function near $0$, $r_k\not=0$.
Let us show that necessarily $i=q, j=1$.

Introduce a new (local) variable $w=L(z)=z-b_k(c)$ and consider
conjugate map
$$g(\rho, w)=L\circ f_c^n\circ L^{-1}(w)$$
where $c=\psi(\rho)$. Then for all $\rho$ near
$\rho_0$ and $w$ near $0$,
\begin{equation}\label{change}
g(\rho,0)=0, \  \frac{\partial g}{\partial w}(\rho, 0)=\rho.
\end{equation}
Now consider $g^q$. Then
$$g^q(\rho_0,w)=w+A w^{q+1}+ O(w^{q+2})$$
where $A\not=0$.
Taking this into account we obtain using (~\ref{change}) that 
\begin{equation}\label{nonzero}
g^q(\rho, w)-w=(\rho^q-1)w+O((\rho-\rho_0)w^2)+Aw^{q+1}+
O(w^{q+2}) \ \ A\not=0.
\end{equation}
On the other hand, 
$b_k(c)-b_k(c_0)=O(\rho-\rho_0)=O(s^i)$ so that a pair of functions
\begin{equation}\label{ij}
(\rho-\rho_0=s^i, w=r_k s^j + O(s^{j+1})+O(s^i))
\end{equation}
is a solution of the equation
$(g^q(\rho, w)-w)/w=0$. 
Substituting this pair of functions of $s$
into (~\ref{nonzero}) we arriv at the conclusion that $i=j q$.

Thus locally (near ($\rho_0$, $b_k(c_0)$))
the algebraic function $(\rho, z)$ has the form
$$\rho-\rho_0=s^q, z-b_k(c_0)=F_k(s)=F_k'(0)s+O(s^2), \ F_k'(0)\not=0$$

Now, $nq$ points $b_k(c_0)+F_k(s\exp(2\pi j/q)), k=1,...,n, j=0,...,q-1$
are fixed by $f_c^{nq}$ and are close to the $n$-cycle
of $f_{c_0}$ with a multiplier a primitive $q$-root of $1$.
Therefore, they form a single $nq$-periodic orbit.
\end{proof}

The following two Lemmas allow us to connect
the multiplier $\rho$ with the multiplier
of a periodic orbit after the bifurcation.
In the first one
we use an idea of~\cite{Che} although our proof is more involved than
in~\cite{Che}. 
The reason is that the limb $L(W, t_0)$ may contain
more than one hyperbolic components of period 
$nq$ or less (which is impossible in case $n=1$ considered in~\cite{Che}).
\begin{lem}\label{q^{-3}} 
There exists $n_0$ as follows. For every
$n>n_0$ and $q>n_0$, if $t_0=p/q$ is $n$-deep,
then each function $F_k$ 
extends to a holomorphic function 
in the disk $\{|s|<9/10\}$.
Moreover, the domain $\psi(B(\rho_0, (9/10)^q))$ 
is disjoint with any limb $L(W, p'/q')$ which is
different from $L(W, t_0)$ and such that $q'\le q+1$.
\end{lem}
\begin{proof}
We prove by several steps that the function $F_k(s)$ extends 
to a analytic function in the disk
$B(0, r^{1/q})$ where $r=1/(n q^{2}(q+1))$.
It will be enough because $(9/10)^q<r$ for big $q>n$

Denote $B=B(\rho_0, r)$ and $\tilde B=\psi(B)$.

1. $\tilde B$ is disjoint with all the limbs $L(W, t')$
of $W$ different from $L(W, t_0)$, such that $t'=p'/q'$ and $q'\le q+1$.
Indeed, let us project the disk $B$ by $\log$.
Since $r$ is very small, the projection is 
(asymptotically) a disk $B(2\pi it_0, r)$.
Doing a simple calculation similar to~\cite{Che} we see that for all $n,q$, 
the latter disk is disjoint 
with all 
Yoccoz's circle that touches the vertical line at a point
$2\pi i p'/q'$ with some $q'\le q+1$.
The claim follows. 

2. Assume the contrary: there exists a path $\gamma$ in 
$B(0, r^{1/q})$ that connect $0$ ans some $s_1$,
such that $F_k$ cannot be extended analytically through $s_1$.
Since $(\rho(s), a(s))$ where $\rho(s)-\rho_0=s^q, z(s)-b_k(c_0)=F_k(s)$
satisfy a polynomial equation, the function $F_k$ has
an analytic continuation along every curve unless it meets a singular point.
Therefore, 
the point $(\rho_1, a)$, such that  $\rho_1:=\rho(s_1)$ and
$a:=b_k(c_0)+F_k(s_1)$ is so that $a$ is a fixed
point of $f^{nq}_{c_1}$  ($c_1=\psi(\rho(s_1))$) 
with the multiplier $1$.

Let us see where $c_1\in B$ can be situated. 
By Lemma~\ref{combin}, $c_1$ belongs either to $L(W, t_0)$
or to some $L(W, p'/q')$ other than $L(W, t_0)$ and such that
$q'\le q+1$. The latter possibility is excluded by Step 1.
Thus $c_1\in L(W, t_0)$.

3. Consider the image $\Gamma$ of the path $\gamma$
by the map $\psi(s^q)$. It connects the root $c_0=c(W, t_0)$
of $L(W, t_0)$ to $c_1$.
The curve $\Gamma$ cannot belong completely to the wake
$W^{**}:=W(t_0)^*$ (which contains the limb $L(W, t_0)$).
The reason is that, as we know, the point
$a(s)$ extends analytically to this wake.
Assume for a moment that there is another curve $\Gamma_1$
that connect $c_1$ and $c_0$ inside the wake $W^{**}$ and such that
it lies in $\tilde B$. Let us deform (kepping the end points fixed)
the curve $\Gamma_1$ to $\Gamma$ inside $\tilde B$. If along the way
we will not meet another singular point, then it contradics to the
fact that the
point $c_1$ is singular for $a(s)$. 
Hence, we must meet another singulat point.
It also must belong to the limb $L(W, t_0)$.
Then we can replave $\Gamma_1$ by another curve that 
connects $c_0$ and $c_1$ inside $W^{**}$ and inside $B$,
such that $\Gamma_1$ can be deformed to $\Gamma$ inside
$B$ without meeting singular points, a contradiction.
 
We conclude that there exists a curve $\Gamma_1$ connecting
$c_0$ and $c_1$ inside the wake $W^{**}$, such that it leaves
$\tilde B$ (and then comes back to $c_1$), such that when we deform
$\Gamma_1$ to $\Gamma$ we meet a singular point $c_2$, which
belongs to another limb $L(W, p_2/q_2)$ with $q_2\le q+1$.
Then the continuum $L(W, p_2/q_2)$ contains the points
$c_2$ and its root and disjoint with $L(W, t_0)$, therefore,
$L(W, p_2/q_2)$ must cross $\tilde B$, a contradiction to the fact
that  $\tilde B$ is disjoint with all such limbs.
\end{proof}

\begin{defi}\label{bifur}
The latter Lemma shows that for every $\rho$ such that
$\rho-\rho_0=s^q$ with $|s|<9/10$, the
points 
$b_k(c_0)+F_k(s\exp(2\pi j/q)), k=1,...,n, j=0,...,q-1$
form a $nq$-cycle of $f_c$ where $c=\psi(\rho_0+s^q)$.
We denote this periodic orbit by $O^q(c)$.
As $c=c_0$, it coincides with the $n$-cycle $O(c_0)$.
\end{defi}

\begin{defi}\label{dist}
If $A, B$ are two sets in the plane, we say
that $A$ is $\delta$-close to $B$ 
and denote it by $d(A, B)<\delta$ iff for every point $a\in A$ there exists
a point $b\in B$, such that $|a-b|<\delta$. 
\end{defi}
\begin{com}
The function $d$ is not symmetric.
On the other hand, it is easy to check the triangle inequality:
if $d(A, B)<\delta_1$ and $d(C, A)<\delta_2$,
then $d(C, B)<\delta_1+\delta_2$.
\end{com}

For $t_0=p/q$, denote
\begin{equation}\label{disk}
B_{t_0}=B(\exp(2\pi it_0), (\frac{9}{10})^q).
\end{equation}

By (~\ref{c}), one can assume that $|c|<3$.
If $z$ is a periodic point of $f_c$ we then have that $|z|\le 6$. 
Now the Schwarz Lemma gives us:
\begin{coro}\label{lip}
For every $\rho\in B_{t_0}$, where $t_0$ is $n$-deep,
and for every $s$ such that $\rho-\rho_0=s^q$ and $|s|< 9/10$,
we have, for $c=\psi(\rho)$,  
we have:
$$d(O^q(c), O(c_0))< 7|s|$$
$$d(O(c), O(c_0))< 6(\frac{10 |s|}{9})^q$$
\end{coro}

\

Denote by $\lambda(\rho)$
the multiplier of the periodic orbit $O^q(c)$ for $c=\psi(\rho)$.
In other words,
$\lambda=\rho_{W(t_0)}\circ \rho^{-1}$ whenever it is well defined.

It follows from, for instance, Lemma~\ref{local} that $\lambda$ is defined
and holomorphic near
$\rho=\rho_0$. Moreover, by~\cite{Gu},
$$\frac{d\lambda}{d\rho}(\rho_0)=-\frac{q^2}{\rho_0}.$$
This formula can be also derived directly from
(~\ref{nonzero}) with help of (~\ref{ij}) where $i=q, j=1$.

As it follows directly from Lemma~\ref{q^{-3}},
for $n$, $q$ biger than $n_0$, the function
$\lambda$ is holomorphic in $B_{t_0}$. But it is not necessarily
univalent there. On the other hand, to choose next internal argument,
we need that the image by $\lambda$ is not small.
This is proved in the next statement.  
\begin{lem}\label{next}
For every real $T\in (-1/2, 1/2]$, the equation 
$\lambda(\rho)=\exp(2\pi i T)$ 
has at most one solution $\rho$ in 
$B_{t_0}$, and for such a solution the
corresponding $c=\psi(\rho)$ lies 
on the boundary of the $nq$-hyperbolic component $W(t_0)$,
which tangents $W$ at the point $c_0=\psi(2\pi it_0)$.
The following covering property
holds:
for every $r\le (9/10)^q$, the image of the disk
$B(\exp(2\pi it_0),r)$ under the map $\lambda$
covers the disk $B(1,q^2 r/16)$.
\end{lem}
\begin{proof}
Let $\lambda(\rho_1)=\exp(2\pi i T_1)$ for some real $T_1$.
Lemma~\ref{q^{-3}} tells us that 
$\psi(B_{t_0})$ is disjoint with any $p'/q'$-limb of $W$
other than $L(W, t_0)$, where
$q'\le q+1$.
On the other hand, $f_{c_1}$ has $nq$-periodic orbit $O^{q}(c_1)$
with the multiplier $\rho_1=\exp(2\pi i T_1)$.
Therefore, $c_1\in L(W, t)$ and moreover lies in the boundary
of an $nq$-hyperbolic component. 
This hyperbolic component belongs to some limb
$L(W, t')$, $t'=p'/q'$, and $q'\ge q+1$. 
If $t'\not=t_0$, then $L(W, t')$ contains a parameter
$\tilde c$ on the boundary of this component, such that
$f_{\tilde c}$ has a $nq$-periodic orbit with the multiplier $1$.
By Lemma~\ref{combin}, $q\ge q'-1$, i.e. $q\ge q'-1$, a contradiction
with the previous inequality $q'\ge q+1$.
Thus $c_1$ is in the limb $L(W, t_0)$.  
Hence, it lies in the boundary of some $nq$-hyperbolic
component belonging this time to $L(W, t_0)$.
On the other hand, the multiplier of the periodic orbit
$O^q(c)$ is bigger on modulus than 1 off the closure
of the component $W(t_0)$.
Thus the hyperbolic component containing $c_1$ in its boundary is just
$W(t_0)$.
Since the corresponding multiplier is
injective on the boundary, it means that $\rho_1$ is the
only solution of the equation $\lambda(\rho)=\exp(2\pi i T_1)$.

To prove the covering property, given $r\le (9/10)^q$,
consider a function
$m(w)=(q^2 r)^{-1} (\lambda(2\pi i t_0)+rw)-1)$.
It is holomorphic in the unit disk, $m(0)=0$,
$|m'(0)|=1$ and, by the proved part of the statement,
$m(w)=0$ if and only if $w=0$.
Therefore, by a classical result (Caratheodory-Fekete, see e.g.~\cite{Gol}),
the disk $B(0, 1/16)$ is covered by the image
of the unit disk under the map $m$.
This is equivalent to the covering property.
\end{proof}

\subsection{Non locally connected Julia sets}
Our aim is to prove Theorem~\ref{intro} stated in the
Introduction. We'll do it in few steps. The main one consists in proving
the following result.
\begin{theo}\label{notlc}
There exists $N$, such that,
for every $n>N$ the following holds.
Let $t_0$, $t_1$,...,$t_m$,... be any sequence
of rational numbers $t_m=p_{m}/q_{m}\in (-1/2, 1/2)$ satisfying
the following property:

(a) if one denotes 
$n_0=n$, $n_m=nq_0q_1...q_{m-1}$ for $m>0$, then
$t_m$ is $n_m$-deep, $m=0,1,2,...$,
see Definition~\ref{deep}.
Moreover, $|t_m|^{1/q_{m-1}}<9/10$, $m=1,2,...$;

Then, for every
$n$-hyperbolic components $W$:

(1) denote $W^0=W$, $W^{m}=W^{m-1}(t_{m-1})$, $m=1,2,...$
(i.e. the hyperbolic component $W^m$ touches the hyperbolic
component $W^{m-1}$ at a point with the internal argument
$t_{m-1}$). Then (a) implies that
the limbs $L(W^{m-1}, t_{m-1})$ form a nested
sequence of compact sets with the intersection a unique 
point $c_*$.

If in addition to (a) we have also:

(b) 
$$\sum_{m=1}^\infty |t_m|^{1/q_{m-1}}< 4^{-n}/16,$$
then:

(2) the map $f_{c^*}$ is infinitely renormalizable
with non locally connected Julia set.
\end{theo}

We begin the proof with some notations.

$$\Delta=4^{-n}/16, \ \ \
\epsilon_m=|t_{m+1}|^{1/q_{m}}, m=0,1,..., \ \ \
n_0=n, \ n_m=nq_0...q_{m-1}, m=1,2,...$$

$$\psi_m(w)=\psi_{W^{m}}(\exp(w)), \ m=0,1,...$$

$$c_m=\psi_m(2\pi i t_m), m=0,1,... .$$
$c_m$ is the point in the boundary of $W^{m}$ with
the internal argument $t_m$.

Since $t_m$ is $n_m$-deep, the function $\psi_m$ 
extends in a univalent fashion to
$B(2\pi it_m, d_m)$.
Remind that $n_m/q_m^4<d_m<4n_m\log 2/q_m$ and
$$B(c_m)=\psi_m(B(2\pi it_m, d_m/2)).$$
By Theorem~\ref{limb}, each $B(c_m)$ is ``roughly''
a round disk around the point $c_m$:
$$B(c_m, r_m/4)\subset B_m\subset B(c_m, 4 r_m)$$
where $r_m=|\psi_m'(2\pi it_m)| d_m/2$.

$$B_m=\psi_m(B(2\pi it_m, (9/10)^{q_m})), \ \
D_m=\psi_m(B(2\pi it_m, |t_{m+1}|)), \ \
D_m'=\psi_m(B(2\pi it_m, 100 |t_{m+1}|))$$

Note that if $r/d_m<1/2$, then
$\psi_m(B((2\pi it_m, r))$ is ``almost'' a disk:
$$B(c_m, 4^{-1} r/r_m)\subset \psi_m(B((2\pi it_m, r))\subset B(c_m, 4 r/r_m)$$
Note also that by the condition (a)
all $q_m$ are big provided $n$ is big.

It implies $100|t_{m+1}|/(9/10)^{q_m}$ 
and $(9/10)^{q_m}/d_m$ are small
(for $n$ big enough),
$\tilde B_m$, $D_m$ and $D_m'$ are ``almost'' disks
(around $c_m$), and for $m=0,1,...$,
$$D_m\subset D_m'\subset B_m\subset B(c_m),$$
and, moreover, the diameter of each previous set is much smaller
than the diameter of the next one.

\begin{lem}\label{incl}
There is $N$, such that, for each $n>N$ the following holds.

(i) $c_{m+1}\in D_m$, $m=0,1...$

(ii) $B(c_{m+1})\subset D_m'$, $m=0,1,...$
\end{lem}
\begin{proof}
We use Lemma~\ref{next}. 
Consider the function $\lambda=\rho_{W^{m+1}}\circ \psi_m$.
Then $\lambda$ is holomorphic
in $B(2\pi it_m, (9/10)^{q_m})$ and 
since $|t_{m+1}|=\epsilon_m^{q_m}<(9/10)^{q_m}$, then
$\lambda(B(2\pi it_m, |t_{m+1}|))$ covers
the disk $B(1, q_m^2 |t_{m+1}|/16)$.
If $q_m>16$ and $t_{m+1}$ is small enough, then
$$\exp(2\pi it_{m+1})\sim 1+2\pi it_{m+1}\in B(1, q_m^2 |t_{m+1}|/16).$$ 
Therefore, 
the point 
$$c_{m+1}=\psi_{m+1}(2\pi i t_{m+1})\in 
\psi_m(B(2\pi it_m, |t_{m+1}|))=D_m$$.

To prove that $B(c_{m+1})$ is contained
in the $100$-times bigger disk $D_m'$,
let us notice that by part (1) of Theorem~\ref{limb},
the root $c_m$ of $W^{m+1}$ is outside of
$B(c_{m+1})$, and
$B(c_{m+1})$ is ``roughly'' a disk around $c_{m+1}$. Hence,
$B(c_{m+1})\subset B(c_{m+1}, 4|c_{m+1}-c_m|)$ and
$B(c_{m+1})\subset B(c_m, 17|c_{m+1}-c_m|)$

Let us estimate $|c_{m+1}-c_m|$.
We use the distortion bounds for univalent maps (~\ref{distor}).
Denote
$\delta_m=2|t_{m+1}|/d_m$. 
Since $c_{m+1}\in D_m$, then 
$|c_{m+1}-c_m|<\delta_m \alpha(\delta_m) r_m$
where $\alpha(x)=(1-x)^{-2}$.
Thus $B(c_{m+1})\subset B(c_m, 17\delta_m \alpha(\delta_m) r_m)$.
On the other hand,
$D_m'=\psi_m(B(2\pi it_m, 100 |t_{m+1}|)$ contains the disk
$B(c_m, 100\delta_m (\alpha(100 \delta_m))^{-1} r_m)$.
We conclude that if $\delta_m<x_0$ where $x_0>0$ is the
solution of equation  
$$100(1-100x_0)^2=17 (1-x_0)^{-2}$$ 
then
$$B(c_{m+1})\subset B(c_m, 17\delta_m \alpha(\delta_m) r_m)
\subset B(c_m, 100\delta_m (\alpha(100 \delta_m))^{-1} r_m)
\subset D_m'.$$

The condition $\delta_m<x_0$ means that
$2|t_{m+1}|< d_m x_0$.
This holds if $(9/10)^{q_m}< (d_m/2) x_0$
which is aparently always the case if
$q_m$ is big enough.
Thus if $n$ is big, the conclusion (ii)
holds. Lemma is proved.
\end{proof}

The lemma implies that for all $m$,
$$D_{m+1}'\subset B_{m+1}\subset B(c_{m+1})\subset D_m'.$$
Since the limb $L(W^m, t_m)$ is contained in $B(c_m)$
and 
the diameters of $D_m'$ is by a definite factor smaller than
the diameter of $B(c_m)$, we conclude
that the limbs $L(W^m, t_m)$
shrink to a point $c_*$. This proves
the conclusion (1) of the Theorem.
Moreover, we have:
$$\{c_*\}=\cap_{m=0}^\infty B_m=\cap_{m=0}^\infty D_m'.$$

Based on this and on Lemma~\ref{lip}, it is not difficult to prove 
the conclusion (2).

Indeed, let us denote by $O_m(c)$ the $n_m$-periodic orbit
of $f_c$, which is attracting if $c\in W^m$, $m=0,1,2,...$.
As we know, $O_m$ extends holomorphically for $c$ in the wake
$(W^m)^*$ of $W^m$.  Moreover, Lemma~\ref{q^{-3}} tells us that
for each $m=0,1,2,...$, the orbit $O_{m+1}$ extends
holomorphically to the disk  $B_m$ near the root $c_{m}$
of the wake $(W^{m+1})^*$.

Since $c_*\subset D_m'\subset B_m$ and $c_{m+1}\subset D_m$, 
Corollary~\ref{lip} allows us to write:
\begin{equation}\label{fin1}
d(O_m(c_*), O_m(c_m))< 6 (\frac{10 (100)^{1/q_m}\epsilon_m}{9})^{q_m}
< 6 (2\Delta)^{q_m}< \Delta
\end{equation}
and
\begin{equation}\label{fin2}
d(O_{m+1}(c_{m+1}), O_m(c_m))< 7 \epsilon_m.
\end{equation}

By the triangle inequality for $d$, we get from (~\ref{fin2}):
$$d(O_{m+1}(c_{m+1}), O_0(c_0))< 7\sum_{k=0}^m \epsilon_k$$
and then (~\ref{fin1}) implies
\begin{equation}\label{away}
d(O_{m+1}(c_*), O_0(c_0))< 7\sum_{k=0}^m \epsilon_k + \Delta
< 8\Delta.
\end{equation}

On the other hand, since $O_0(c_0)$ is a neutral $n$-periodic orbit
of $f_{c_0}$, the distance of $O_0(c_0)$ to zero is bigger than
$4^{-n}$. Hence, we have
if $\Delta=4^{-n}/16$, then, for all $m=0,1,2,...$,
every point of $O_m(c_*)$ is away from the origin by the distance
at least $4^{-n}/2$.

It is known that this implies the
non local connectivity of $J_{c_*}$ (see~\cite{Sor} for a detailed proof,
although our description is closer to~\cite{Mil}).
Let us outline a proof. 
Denote by $\tau_m, \tau_m'$ two external arguments
in the $c$-plane to the root $c_m$ of the wake of $W^m$.
By the formula (~\ref{digits}),
$|\tau_{m+1}-\tau_{m+1}'|\le (2^{n_m})^2/(2^{n_m q_m}-1)\to 0$
as $m\to \infty$, that is $\tau_m, \tau_m'$ tend to some $\tau^*$.

On the other hand, 
for every $c$ in the wake of $W^{m+1}$, we have the
following picture in the dynamical plane
of $f_c$ (see beginning of the present Section): 
two external rays of $f_c$ with arguments
$\tau_{m+1}$ and $\tau_{m+1}'$ land at a single point $b_m(c)$, which
is a point of the repelling periodic orbit $O_m(c)$.
Moreover,
$c$ lies in the sector 
bounded by these two external rays 
and disjoint with the origin.

Now assume that the Julia set of $f_{c_*}$ is locally connected.
Then this discussion implies that, in the dynamical plane of $f_{c_*}$, 
the external argument
of $c_*$ is equal to $\tau_*$ and since $\tau_m\to \tau_*$, also
$b_m(c_*)\to c_*$. Taking a preimage by $f_{c_*}$, it gives us 
a sequence of points of the sets $O_m(c_*)$ which tends 
to the origin, a contradiction with the fact that the orbits
$O_m(c_*)$ stay away from the origin.
\subsection{Proof of Theorem~\ref{intro}}
Here we prove a more general Theorem~\ref{introcopy}.
Theorem~\ref{intro} announced in the Introduction
is an immediate corollary of this together
with Definition~\ref{deep} and sufficient conditions
(~\ref{first})-(~\ref{second}).
\begin{theo}\label{introcopy}
Let $n\ge 1$. Let
$$t_0, t_1,..., t_m,...$$
be any sequence
of rational numbers $t_m=p_m/q_m\in (-1/2, 1/2]$ which satisfy
the following properties:
if we denote $n_0=n$, $n_m=nq_0...q_{m-1}$, $m>0$, then,
for all $m$ large enough,
$t_m$ is $n_m$-deep,
and also
\begin{equation}\label{assump2copy}
\sum_{m=1}^\infty |t_m|^{1/q_{m-1}}<\infty.
\end{equation}
Given a hyperbolic component $W$ of the Mandelbrot set
of period $n$,
consider a sequence of hyperbolic
components $W^m$:  $W_0=W$, and, for $m>0$,
$W^m=W^{m-1}(t_{m-1})$, i.e. $W^m$ touches the hyperbolic
component $W^{m-1}$ at a point with the internal argument
$t_{m-1}$.
For every $m$, consider the $t_m$-limb $L(W^m, t_m)$ of $W^m$
(it contains $W^{m+1}$).
Then the sequence of limbs $L(W^m, t_m)$ shrink 
to a unique point $c_*$,
the Mandelbrot set is locally connected at $c_*$, and
the map $f_{c_*}$ is infinitely renormalizable
with non locally connected Julia set.
\end{theo}

Let $m_0$ be large enough, so that $|t_{m+1}|<(9/10)^{q_m}$
and $t_m$ is $n_m$-deep for every $m\ge m_0$.
Then the sequence $t_{m_0},...,t_m,...$ satisfies
the conditon (a) of Theorem~\ref{notlc} with $n=n_{m_0}$. 
Hence, by the conclusion (1) of Theorem~\ref{notlc}
the limbs $L(W^m, t_m)$ shrink to a unique point $c_*$.

It remains to show that the Julia set $J_{c_*}$
of $f_{c_*}$ is not locally connected. 
We show that {\it whether $J_{c_*}$ is locally connected or not
depends only on a tail of the sequence $t_0,...,t_m,...$.}
Let us introduce the following
notation. Let $W_0$ be the $1$-hyperbolic component (the main cardioid).
Given $k$, let's start with the component $W_0$ and the
tail $T_k=\{t_k, t_{k+1},...\}$ in place of $W$ and $\{t_0,t_1,...\}$.
Then we get a sequence of hyperbolic components
$W_0^{k,m}$, $m\ge k$, where $W_0^{k,k}=W_0$ and,
for $m>k$, the component $W_0^{k,m}$
touches $W_0^{k,m-1}$ at the point with internal argument $t_{m-1}$.
We have proved that, for every $k\ge m_0$, the sequence of limbs
$L(W_0^{k,m})$, $m=k,k+1,...$ shrinks to
a unique point $c_*^k$. 
Now, we have
\begin{prop}\label{mcard}
For every $k$ large enough,
the Julia set of $f_{c_*^k}$ is not locally connected.
\end{prop}
\begin{proof}
If $c$ is in the boundary of the main cardioid, the 
neutral fixed point $\alpha(c)$ of $f_c$ has modulus $1/2$.
Since $t_k\to 0$,
by continuity (and Yoccoz's bound), there is $\delta_0>0$,
such that, for every $k$ large enough,  
for $c$ in the boundary of $W(t_k)$, the 
neutral periodic orbit of $f_c$ (of period $q_k$)
stays away from $0$ by at least $\delta_0$.
If we increase $k$ further, 
one can assume that
$$\sum_{m=k}^\infty |t_{m+1}|^{1/q_m}<4^{-q_k}/16$$
and also
$4^{-q_k}/2<\delta_0$. Define $\delta=\delta_0-4^{-q_k}/2$.
Now 
we can apply Theorem~\ref{notlc} to $q_k$ instead of $n$ and to the sequence
$t_{k+1}, t_{k+2},...$
taking into account (~\ref{away}). It shows that
all corresponding periodic orbits
at the parameter $c_*^k$ will stay away from
the origin by at least $\delta$.  
\end{proof}
Let us come back to the map $f_{c_*}$.
Next statement will finish the proof that $J_{c_*}$ is not
locally connected.
\begin{prop}\label{tuning}
For every $k$, there is a restriction $f_{c_*}^{n_k}: U_k\to V_k$,
such that this is a polynomial-like map of degree $2$, which
is quasi-conformally conjugate to $f_{c_*^k}$.
\end{prop}
\begin{proof}
Consider first the wake of $W^{k-1}$. For $c$ in this wake,
$f_c$ has a periodic orbit $O_{k-1}(c)$ of period $n_{k-1}$
which is attracting iff $c\in W^{k-1}$ and is holomorphic
in $c$.   
Now consider the wake $(W^k)^*$ of the next component
$W^k$. Denote also $a_k$ the root of $W^k$.
If $c=a_k$, then $O_{k-1}(a_k)$ is neutral with the 
multiplier $\exp(2\pi i t_{k-1})$: there are $q_{k-1}$
petals attached to each point of  $O_{k-1}(a_k)$, and
$f_{a_k}^{n_{k-1}}$ acts (locally) on them as a rotation 
with the totation number $t_{k-1}$.
If $c\in (W^k)^*$, then  
$O_{k-1}(c)$ is repelling, and has the rotation number $t_{k-1}$.
There exists a point $b_{k-1}(c)$ of $O_{k-1}(c)$, such that
for $c\in   (W^k)^*\cup \{a_k\}$, there are two external rays $R_k(c)$,
$\tilde R_k(c)$ of arguments
$\tau(a_k)$, $\tilde\tau(a_k)$ landing
at $b_{k-1}(c)$, such that the component of the plane
bounded by these two rays that contains $c$ contains no
other ray to the orbit. There are two topological disks
$c\in S'_k(c)\subset S_k(c)$, such that, $S_k(c)$ is bounded by $b_{k-1}(c)$, 
the rays $R_k(c)$, $\tilde R_k(c)$ and by some equipotential,
and $f_c^{n_k}: S'_k(c)\to S_k(c)$ is  a proper map of degree $2$.
By a ``thickening'' of $S_k(c)$, one can turn it into
a polynomial-like map $P_{k, c}: U_{k,c}\to V_{k,c}$,
see~\cite{Mi2} for details.

If $c\in B:=\{c_*\}\cup \cup_{m=k}^\infty (W^m\cup \{a_k\})$, 
we claim that the 
Julia set of $P_{k, c}: U_{k,c}\to V_{k,c}$ is connected.
Indeed, consider the iterates $f_c^{n_k j}(c)$, $j=0,1,...$.
If $c$ is close to $a_k$ and in $W^k$, then
they converge to $b_k(c)$. On the otehr hand, neither of these
iterates cannot cross the boundary of $S'_k(c)$ when
$c\in B$, because $f_c$ is not
Misiurewich map for such $c$.
This proves the claim.
 
By the Straightening Theorem~\cite{DH3},
for every $c\in B$, there is a unique $\nu(c)\in M$, such that
$P_{k,c}: U_{k,c}\to V_{k,c}$ is hybrid equivalent to $f_{\nu(c)}$.
Moreover, by~\cite{DH3}, the map $c\mapsto \nu(c)$ is continuous
(it follows essentially
from the compactness of $K$-quasiconformal maps and
Proposition 7 of~\cite{DH3}).
We need to show that $\nu(c_*)=c_*^k$.
By continuity, it is enough to show that, for every root $a_m$ with $m>k$,
$\nu(a_m)$ is the root $a_m^k$ of the hyperbolic component
$W_0^{k,m}$. 
Let us prove it by induction in $m$.

(1) $m=k+1$. Notice that if $c$ is close to $a_k$ and in $W^{k}$, then
$P_{k, c}: U_{k,c}\to V_{k,c}$ has an attracting fixed point: the
point $b_k(c)$ of $O_{k}(c)$ which lies in $S'_k(c)$ and coincides with
$b_{k-1}(c)$ when $c=a_k$. This attracting fixed point persists as
$c\in W^k$ (periodic points cannot leave $S'_k(c)$), and when
$c=a_{k+1}$, it has the multiplier $\exp(2\pi it_{k})$. Hence,
the conjugate map $f_{\nu(a_{k+1})}$ has a neutral fixed point, and
it interchanges the petals at this point with the same
rotation number $t_{k}$. It follows,
$\nu(a_{k+1})=a_{k+1}^k$, the unique point in the boundary
of the main cardioid with the internal argument $t_k$.

(2) assume that $\nu(a_m)=a_m^k$ holds for some $m\ge k+1$. 
Then a similar argument shows that $\nu(a_{m+1})=a_{m+1}^k$. 
\end{proof}
\subsection{Maps with unbounded combinatorics}
Let us consider the case $n=1$, that is we start with
the main cardioid. Then we can reformulate Theorem~\ref{intro}
as follows.
\begin{theo}\label{infi}
Let a sequence $p_m/q_m$, $m=0,1,...$, of non-zero
rational numbers in $(-1/2, 1/2)$ be such that:
$$\sum_{m=1}^\infty |\frac{p_m}{q_m}|^{1/q_{m-1}}< \infty,$$
$$p_m > (2B_0)  4^{n_m}, \ \ \ \ 2n_m^2 < \frac{p_m^2}{q_m},$$
for all $m$ large enough, where $n_m=q_0...q_{m-1}$.

Let $f_c$ be an infinitely renormalizable polynomial with
the following combinatorial data.

(1) The renormalization periods consists of the sequence 
$\{n_m\}_{m=1}^\infty$.

(2) Denote $J_0=J(f_c)$, and, 
for every $m>0$, $J_m$ be the Julia set
of the renormalization of period $n_m$, which contains $0$.
Let $\alpha_m$ and  
$\beta_m$ be the $\alpha$ (i.e. separating) and
$\beta$ (i.e. non separating)
fixed points of $f_c^{n_m}: J_m\to J_m$. 
We assume that, for every $m$, $\beta_{m+1}=\alpha_{m}$,
and the rotation number of 
$\alpha_m$ w.r.t. the map $f_c^{n_m}: J_m\to J_m$ is equal to $p_m/q_m$.

Then there exists a unique polynomial which satisfies the conditions
(1)-(2), and its Julia set is not locally connected.  
\end{theo}

\end{document}